\documentclass[a4paper,12pt]{article}
\usepackage[dvips]{epsfig}
\usepackage{amsmath,amssymb,amsbsy, amstext,amscd,amsfonts}
\usepackage{color,enumerate,euscript,graphicx,hyperref,srcltx,tikz,indentfirst}
\usepackage{tikz}
\usepackage{xcolor}
\usetikzlibrary{shapes.geometric}
\usetikzlibrary{ positioning,  shapes.geometric}

\input amssym.def
\input amssym.tex

\newtheorem{lemma}{Lemma}[section]%
\newtheorem{theorem}[lemma]{Theorem}%
\newtheorem{proposition}[lemma]{Proposition}%
\newtheorem{definition}[lemma]{Definition}%

 \def\op{\overline P} 
  
   \def\o1{\overline 1}

 \def\ola{\overline a} 
 \def\olz{\overline z}  \def\olb{\overline b}

\def\lg{\langle} \def\rg{\rangle} 
  \def\PGL{\hbox{\rm PGL}}\def\PG{\hbox{\rm PG}}
\def\Aut{\hbox{\rm Aut\,}} 
 \def\cal{\mathcal} 
\def\Cay{\hbox{\rm Cay}}

\def\ZZ{\mathbb{Z}}\def\FF{\mathbb{F}}\def\f{\noindent}
\def\b{\beta} \def\di{\bigm|}  
\def\s{\sigma}\def\o{{\rm o}}\def\AGL{\hbox{\rm AGL}} 
  \def\th{\theta}
\def\pf{\noindent{\it Proof.} } \def\Syl{\hbox{\rm Syl}}
\def\qed{\hfill $\Box$} \def\demo{\pf}
 
\def\GL{\mathrm{GL}} \def\AGL{\mathrm{AGL}}

\def\rtimes{:}

\def\si{\Sigma}

\begin{document}
\begin{center}
{\bf\large Hamilton cycles of semisymmetric graphs of order $2p^3$}
\footnote{This work was supported by the National Natural Science Foundation of China(12301446).}
\end{center}
\begin{center}
Huye Chen\\

{\it {\small School of Mathematics,
Guangxi University,
Nanning 530004, P. R. China.}}\\
{\it {\small Email: chenhy280@gxu.edu.cn}}
\end{center}

\renewcommand{\thefootnote}{\empty}


 \begin{abstract} In light of Lov\'{a}sz's longstanding question on the existence of Hamilton paths in vertex-transitive graphs, Du and Yuan \cite{DY} considered a natural variant: what if vertex-transitivity is relaxed, while a high degree of symmetry--specifically edge-transitivity--is retained? To investigate this, they studied  semisymmetric graphs (i.e. regular, edge-transitive, but not vertex-transitive graphs) and showed that every connected semisymmetric graph of order $2pq$, where $p$ and $q$ are distinct primes, contains a Hamilton cycle. In this paper, it is shown  that for any prime $p$, every connected semisymmetric graph of order $2p^3$
also contains a Hamilton cycle.
  \end{abstract}

\section{Introduction}\label{Introduction}
In this introductory section, we briefly survey the current research status of Lov\'{a}sz's long-standing question on Hamilton paths in vertex-transitive graphs and that of  semisymmetric graphs;
and  then present the main theorem of this paper.
\subsection{Lov\'{a}sz's Question}
A {\it Hamilton path} (resp. {\it cycle}) is a simple path (resp. cycle) that visits every vertex of the graph exactly once. A simple graph containing a Hamilton cycle is called a {\it Hamiltonian graph,}
  In 1969, Lov\'asz \cite{L70} posed the following question:
   \vskip 2mm {\it Do there exist connected  vertex-transitive graphs with no Hamilton paths?}\
   \vskip 2mm
    \f Later, in 1981, Alspach \cite{A81} further asked  whether there exist infinitely many connected vertex-transitive graphs without Hamilton cycles.
 These two questions bridge two seemingly unrelated concepts: the traversability and symmetry of graphs.
To date, there are only four connected vertex-transitive graphs (with at least three vertices) that contain no Hamilton cycle:
the Petersen graph, the Coxeter graph and the truncation of these two graphs.

It has been shown that every connected  vertex-transitive graph (abbreviated as CVTG) of prime order $p$ contains a Hamilton cycle, see \cite{A81}.
For CVTGs of order $n=pq$ (where $p$ and $q$ are primes), partial results were established by Alspach \cite{A79,A81} for $n\in\{2p,3p\}$, and later extended to $n=5p$ in \cite{MP82}. It was eventually shown by Du, Kutnar and Maru\v si\v c that every CVTG of such order contains a Hamilton cycle, except for the Petersen graph, see \cite{DuH-pq}.
Partial results also exist for CVTGs (or vertex-transitive digraphs) of order $n=2pq$ and $p^k$ for small $k$ (where $p$ and $q$ are primes); see \cite{C98,DM83,MP82,Z15}, the survey paper \cite{KM09}, and references therein.
With the obvious exception of the Petersen graph, Hamilton cycles are also known to exist in connected vertex-transitive graphs whose automorphism groups contain a transitive subgroup with a cyclic commutator subgroup of prime-power order \cite{DGMW}.
It is known that Cayley digraphs on $p^k$ vertices all contain Hamilton cycles; see \cite{Wit}. For survey papers on Hamilton cycles in Cayley graphs, see \cite{KM12,Wit1}.

\vskip 3mm
Motivated by Lov\'{a}sz's long-standing question on the existence of Hamilton paths in vertex-transitive graphs, Du and Yuan \cite{DY} considered a natural variant: what if vertex-transitivity is relaxed, while a high degree of symmetry--specifically edge-transitivity--is retained? To investigate this, they studied {\it semisymmetric} graphs (regular, edge-transitive, but not vertex-transitive graphs), noting that every such graph is bipartite with two equal parts (see next subsection), and showed that every connected semisymmetric graph of order $2pq$ (where $p$ and $q$ are distinct primes) contains a Hamilton cycle. Since no semisymmetric graphs exist of order $2p$ or $2p^2$, a natural candidate for this problem is the class of connected semisymmetric graphs of order $2p^3$, where   $p$ is a  prime. In this paper, we prove that every such graph contains a Hamilton cycle.

\subsection{Semisymmretic graphs}
A graph is said to be {\em semisymmetric} if it is regular and
edge-transitive but not vertex-transitive.
It is easy to see that every semisymmetric graph is bipartite, with the two parts being equal in size and each part being vertex-transitive under the edge-transitive group. The smallest semisymmetric graphs have order $20$.

The first person who studied semisymmetric graphs was Folkman. In 1967, he constructed several infinite families of such
graphs and proposed eight open problems, see \cite{Fol}.
Since then, much work has been done on semisymmetric graphs, see \cite{Bou2,II,Kli}. These works provided new constructions of such graphs and resolved nearly all of Folkman's open problems.
In particular, using group-theoretical methods, Iofinova and Ivanov \cite{II} classified in 1985 the cubic
semisymmetric  graphs whose automorphism group acts primitively on both biparts.
This was the first classification
result for such graphs.
 More recently, following some deep results in group theory that rely on the classification of
finite simple groups, together with certain methods from graph coverings, several new  results on semisymmetric graphs have appeared, see
\cite{CMMP,DX} for example.

Folkman \cite{Fol} proved that there are no semisymmetric graphs of order $2p$ or $2p^2$, where $p$ is a prime. Du and Xu \cite{DX} classified semisymmetric graphs of order $2pq$ for two distinct primes $p$ and $q$.
Therefore, a natural question is to determine the semisymmetric graphs of order $2p^3$ (with $p$ prime).
 It was proved in \cite{MMW} that the Gray graph of order $54$ is the only cubic semisymmetric graph of order $2p^3$.
  The classification of all semisymmetric graphs of order $2p^3$ remains an attractive and difficult open problem. These graphs $X$ are naturally divided into two subclasses:
 $\Aut(X)$ acts faithfully on both  biparts; and  unfaithfully on at least one bipart.
The second subclass have been classified by Du and Wang, see \cite{Du-Wang}.
Their classification  will be used in this paper. Since the first subclass remains open, we need to prove the existence of Hamilton cycles by a direct approach.

Now we are ready to state our main theorem in this paper.
\begin{theorem}\label{main}
For any prime  $p$, every connected semisymmtric  graph of order $2p^3$  contains a Hamilton cycle.
\end{theorem}

\section{Preliminaries}
Throughout this paper, all graphs are finite, undirected, and simple unless explicitly indicated. Henceforth, $p$ denotes an odd prime, recalling that   the smallest semisymmetric graphs have order $20$.
Given a graph $X$, by $V(X),~E(X)$ and $\Aut(X)$ we denote the vertex set, the edge set and the automorphism group of $X$, respectively. If $X$ is bipartite with bipartition $V(X)=W(X)\cup U(X)$, we let $\Aut(X)^+:=\lg g\in\Aut(X)\mid U(X)^g=U(X), W(X)^g=W(X)\rg.$
For disjoint subsets $U$ and $W$ of $V(X)$, $X[U,W]$ denotes the induced bipartite subgraph with bipartition $U$ and $W$.
In the case when $X[U,W]$ is regular, $d(U,W)$ denotes the valency of $X[U,W]$.

Let $G$ be a permutation group on a set $V$. For any subset $A\subseteq V$,
$G_A$ (resp. $G_{(A)}$) denotes the set-wise (resp. point-wise) stabilizer of $A$ in $G$.
A {\it block} (or {\it $|B|$-block}) of a transitive permutation group $G$ is a nonempty subset $B\subseteq V$ such that for every $g\in G$, either $B^g=B$ or $B^g\cap B=\emptyset$.
For a block $B\subseteq V$, the collection $\mathcal{B}=\{B^g\mid g\in G\}$ forms a {\it block system} of $G$ on $V$. Obviously, given a   normal    subgroup $N$ of $G$,
 all  $N$-orbits form  a block system of $G$.

\vskip 3mm
\f {\bf 2.1 Semiregular automorphisms and  lifting cycle technique}

\vskip 3mm
Let $m\geqslant 1$ and $n\geqslant 2$ be integers. An automorphism $\rho$ of a graph $X$ is called $(m,n)$-{\em semiregular} (in short, {\em semiregular}) if, as a permutation on $V(X)$, its cycle decomposition consists of exactly $m$ cycles, each of length $n$.
Let $\mathcal{P}$ be the set of orbits of $\lg \rho \rg $, where $\rho$ is a semiregular automorphism.
Let $X_{\mathcal{P}}$ be the {\em quotient graph} corresponding to $\mathcal{P}$, the graph  whose vertex set is $\mathcal{P}$, with $A, B \in \mathcal{P}$ adjacent if there exist adjacent vertices $a \in A$ and $b \in B$ in $X$.
Let $X_\rho$ be the {\em quotient multigraph} corresponding to $\rho$,  the multigraph  whose vertex set is $\mathcal{P}$ and in which $A,B \in \mathcal{P}$ are joined by $d(A,B)$ edges.
Note that the quotient graph $X_\mathcal{P}$ is precisely the underlying graph of $X_\rho$.

Let $X$ be a graph with $G\leq\Aut(X)$, and $N\unlhd G$ a nontrivial normal subgroup. Let $\mathcal{B}$ be the set of $N$-orbits on $V(X)$. The {\em normal quotient graph} $X_{N}$ of $X$ induced by $N$ is defined with vertex set $\mathcal{B}$, and two vertices $B,B'\in \mathcal{B}$ are adjacent if and only if there exist vertices $\b\in B$ and $\b'\in B'$ that are adjacent in the graph $X$.

A key tool for addressing Hamilton cycle problems is the lifting cycle technique (see \cite{A1989,KM09,DM83}).
When the quotient graph is applied relative to a semiregular automorphism of prime order and the corresponding quotient multigraph possesses two adjacent orbits linked by a double edge encompassed within a Hamilton cycle, lifts of Hamilton cycles from quotient graphs are invariably achievable. This double edge enables us to conveniently ``change direction" to procure a walk in the quotient graph that elevates into a full cycle above.

Let $X$ be a graph admitting a $(m,n)$-semiregular automorphism $\rho$.
Let $\mathcal{P} = \{S_1, S_2, \cdots , S_m\}$ be the set of orbits of $\langle\rho\rangle$, and let $\pi : X \to X_{\mathcal{P}}$ be the corresponding projection of $X$ to its quotient graph $X_{\mathcal{P}}$. For a (possibly closed) path $W = S_{i_1}S_{i_2}\ldots S_{i_k}$ in $X_{\mathcal{P}}$ we let the {\em lift} of $W$ be the set of all paths in $X$ that project to $W$. The proof of following lemma is straightforward and is just a reformulation of \cite[Lemma~5]{MP82}.

\begin{proposition} \label{pro:4}
Let $X$ be a graph admitting
a $(m,p)$-semiregular automorphism $\rho$, where $p$ is a prime.
Let $C$ be a cycle of length $k$ in the quotient graph $X_{\mathcal{P}}$,
where $\mathcal{P}$ is the set of orbits of $\langle\rho\rangle$.
Then, the lift of $C$ either contains a cycle of length
$kp$ or it consists of $p$ disjoint $k$-cycles.
In the latter case we have $d(S,S') = 1$ for every edge $SS'$ of $C$.
\end{proposition}

\vskip 3mm
\f {\bf 2.2 Some results  on graph theory and group theory}
\vskip 3mm


 If $X$ is bipartite with its two parts corresponding to two $G$-orbits, then $X$ is isomorphic to a bi-coset graph of $G$.
We now formally define  bi-coset graphs and bi-Cayley graphs.

\begin{definition}\label{coset-graph}
Let $G$ be a group with two subgroups $L$ and $R$, while by    $[G:L]$ and  $[G:R]$  we denote the set of right cosets of $G$ relative to $L$ and $R$, respectively.
 Let $D$ be a union of some double cosets of $L$ and $R$ in $G$, namely, $D=\cup_iRd_iL$.
The bi-coset graph ${\bf B}(G;L,R,D)$ is the bipartite graph with
       vertex partition: $[G:L]\cup [G:R]$  and the  edge set: $\{(Lg,Rdg)\mid g\in G, d\in D\}$.

     Moreover, if $L=R=1$, then $D$ is a subset of $G$ and $X$ is called a bi-Cayley graph of $G$ with respect to $D$.
\end{definition}

\begin{proposition}{\rm\cite[Lemmas 2.3, 2.5,2.6]{DX}}\label{coset-graph-pro}
Using the notations from Definition \ref{coset-graph}, we have the following properties$:$
\begin{enumerate}
  \item The bi-coset graph $X:={\bf B}(G;L,R,D)$ is connected if and only if $\lg D^{-1}D \rg=G.$
  \item The graph $X$ is $G$-edge-transitive if and only if $D=RdL$ for some $d\in G.$
  \item A bi-coset graph $X$ is bi-Cayley if and only if $\Aut(X)^+$ has a subgroup which acts regularly on each of the two parts of $V(X).$
  \item If $G$ is abelian and acts regularly on both parts of $X$, then $X$ is vertex-transitive. In other words, bi-Cayley graph of abelian groups are vertex-transitive.
\end{enumerate}
\end{proposition}


\begin{proposition}{\rm\cite{Wit}}\label{cayleyp3}
Every connected Cayley digraph on a group of prime-power order greater than $2$ has a Hamilton cycle.
\end{proposition}

\begin{proposition}{\rm\cite{A79,DM87}}\label{2p^2DM}
Every connected vertex-transitive graph of order $2p$ or $2p^2$ has a Hamilton cycle.
\end{proposition}



\begin{proposition}{\rm\cite[I. Satz.7.2]{Hup}}\label{pcentral}
If $G$ is a finite $p$-group and $1\neq H\unlhd G$, then $H\cap Z(G)\neq 1$.
\end{proposition}

\section{Proof of Theorem~\ref{main}}
\demo Let $X$ be a connected semisymmetric graph with automorphism group $A$. Let $V=W\cup U$ be a partition of vertex set $V=V(X)$.
Let $P\in \Syl_p(A)$, where $\Syl_p(A)$ denote the set of all Sylow $p$-subgroups of $A$.
Then we shall  discuss two cases, separately: $A$ acts faithfully on both  biparts  in Section 4 and unfaithfully on one bipart in Section 5.
It will be shown that for both cases, $X$ has a Hamilton cycle, see Theorems~\ref{main1} and ~\ref{main2}.  Thus, Theorem~\ref{main} is proved. \qed

\vskip 3mm
To end up this section,  two general lemmas are presented.

\begin{lemma}\label{main-lemma}
Let $X$ be a graph and let $G\leq \Aut(X)$ be a subgroup. Let $\mathcal{B}$ be a  set of orbits of an $(m,p)$-semiregular automorphism $\rho\in \Aut(X)$, where $p$ is a prime and $m$ is a positive integer. Suppose $\mathcal{B}$ is invariant under the action of $G$. Then $X$ contains a Hamilton cycle, provided that the quotient graph $X_{\mathcal{B}}$ admits a Hamilton cycle and that $K_{(B_i)}\neq 1$ for some $i\in\ZZ_m$, where $K=G_{(\mathcal{B})}$ is the kernel of the action of $G$ on $\mathcal{B}$.
\end{lemma}
\demo
Suppose that $\mathcal{C}$ is a Hamilton cycle of the quotient graph $X_{\mathcal{B}}$,
where $$\mathcal{C}:B_0\sim B_1\sim \cdots\sim B_{m-1}\sim B_0,$$
with $|\mathcal{B}|=m$, $B_i\in\mathcal{B}$, and
$i\in\ZZ_m.$
Then
 each induced bipartite subgraph $X[B_j,B_{j+1}]$ is either a perfect matching or not, where $j\in\ZZ_m$.
If there exists $\jmath\in\ZZ_m$ such that $X[B_\jmath,B_{\jmath+1}]$ is not a perfect matching, then the graph $X$ admits a Hamilton cycle by Proposition \ref{pro:4}.
If each induced subgraph $X[B_\jmath,B_{\jmath+1}]$ is a perfect matching with $\jmath\in\ZZ_m$, then consider the action of $K_{(B_\imath)}\neq 1$ on the sets $B_j\in\mathcal{B}$ ($j\in\ZZ_m$) along the Hamilton cycle $\mathcal{C}$ staring at the set $B_{\imath}$. This action forces $K_{(B_\imath)}=1$, a contradiction.
\qed

\begin{lemma}\label{p3}
 Let $X$ be a semisymmetric graph of order $2p^3$.  Suppose that the automorphism group $A=\Aut(X)$  contains a subgroup $P_1$ that acts faithfully and regularly on both bipartite parts. Then $X$ has a Hamilton cycle.
\end{lemma}
\demo
Suppose that $A$ contains a subgroup $P_1$ acting faithfully and regularly on both bipartite parts.
 If $P_1$ is abelian, then  by Proposition~\ref{coset-graph-pro}, $X$ is vertex-transitive, a contradiction (since $X$ is semisymmetric).
Therefore, we assume that $P_1$ is nonabelian and $|P_1|=p^3$.
Then $X$ can be viewed as a bi-Cayley graph over the nonabelian group $P_1$ with some connection set $S$. Without loss of generality, we assume that $S=\{1,S_1'\}$ and that $P_1= \langle S_1' \rangle$, where $1\notin S_1'\subseteq P_1$. Hence $X$ is a bi-Cayley graph $X\cong \operatorname{B}(P_1, 1, 1, S)$.

Set $V(X)=W(X)\cup U(X)$  with $W(X)=\{W_{x}\mid x\in P_1\}$ and $U(X)=\{U_{x}\mid x\in P_1\}$.
Now consider the Cayley digraph $X'=\Cay(P, S_1')$.
By Proposition~\ref{cayleyp3}, there exists a Hamilton di-cycle $\overrightarrow{C}$ in the digraph $X'$,
where $\overrightarrow{C}: x_0\rightarrow x_1\rightarrow x_2\rightarrow\cdots\rightarrow x_{p^3-1}\rightarrow x_0,$ with $P_1=\{x_i\mid i\in\ZZ_{p^3}\}$.
This directed cycle will be used to construct a Hamilton cycle $\mathcal{C}$ in the graph $X$:
$$\hskip 1.7cm \mathcal{C}: W_{x_0}\sim U_{x_1}\sim W_{x_1}\sim U_{x_2}\sim W_{x_2}\sim\cdots\sim W_{x_{p^3-1}}\sim U_{x_{0}}\sim W_{x_0}. \hskip 1.3cm \Box$$

\section{$A$ acts faithfully on both biparts}
Now, we suppose that $A$ acts faithfully on both bipartite parts in this section.
Then $P\in\Syl_p(A)$ acts faithfully on both bipartite parts too.
 Let $Z=Z(P)$ with $|Z|=p^i$, $i\in\{1,2,3\}$.

 If $i=3$, then $Z$ is regular on both biparts,  forcing $X$ is vertex-transitive  by  Proposition~\ref{coset-graph-pro}, a contradiction.
  Suppose that $i=2$. Then we  consider the quotient graph $X_{ Z}$ induced by $Z$-orbtis on two biparts,  where $V(X_{ Z})=\mathcal{W}\cup\mathcal{U}$. Consider the action of $P$ on
  $X_{ Z}$, whit the kernel $P_0$. Then $P/P_0\cong \ZZ_p$ or $\ZZ_P^2$. In both cases,   there exist an  element $\rho P_0$ of $P/P_0$ moving  $Z$-orbits  on both  two  biparts, so that  the abelian subgroup $\lg Z, \rho\rg$ acts transitively   on both $W$ and $U$.
  Since $A$ acts faithfully on both bipartite parts, we know that   $\lg Z, \rho\rg$ acts  regularly   on both $W$ and $U$,  a contradiction agian.
    In summary,   we get  $|Z|=p$. In particular, we may further assume that $P$ is nonabelian.

Let ${\cal B}_w$  and ${\cal B}_u$ be the set of $Z$-orbits on the two biparts $W$ and $U$, respectively. Denote by $Y:=X_{Z}$ the  quotient graph of $X$ induced by $Z$-orbits. For the action of $P$ on ${\cal B}_w$ and ${\cal B}_u$, let $K_w$ (resp. $K_u$) be the kernel of the action on ${\cal B}_w$ (resp. ${\cal B}_u$). Define $K := K_w \cap K_u$. Then $Z \leq K$.
If $|P|=p^3$, then by Lemma~\ref{p3}, $X$ is Hamiltonian. From now on, we assume $|P|\ge p^4$.

The main result of this section is the following theorem.
\begin{theorem}\label{main1}
 Using the above notation,  suppose that $A$ acts faithfully on both biparts. Then $X$ contains a Hamilton cycle.
\end{theorem}
From now on, let $\op=P/K.$
We shall deal with the following cases: $\op$ acts faithfully on both ${\cal B_w}$ and ${\cal B_u}$ (i.e., $K=K_w=K_u$) in Section 4.1;
$\op$ acts unfaithfully on both ${\cal B_w}$ and ${\cal B_u}$ (i.e., $K\lneq K_w$ and $K\lneq K_u$) in Section 4.2; and $\op$ acts unfaithfully on exactly one of ${\cal B_w}$ and ${\cal B_u}$
(i.e., either $K_w\lneq K_u$ or $K_u\lneq K_w$) in Section 4.3. Therefore, Theorem~\ref{main1} will be proved by using Lemmas~\ref{main11}, \ref{main12} and \ref{main13}.

\subsection{$\op$ acts faithfully on both ${\cal B_w}$ and ${\cal B_u}$}
 We first present a group-theoretic lemma.

\begin{lemma}\label{twofaithful}
Let $Q$ be a $p$-group  having two faithful transitive representations of degree $p^2$ on $W'$ and $U'$, respectively.
Then $Q$ must contain a subgroup  $Q_1\cong\ZZ_p^2$ or  $\ZZ_{p^2}$,   which acts regularly on both $W'$ and $U'$.
\end{lemma}
\demo   Take $z_0'\in Z(Q)$ such that $o(z_0')=p$. Then $\lg z_0'\rg$  induces  two complete $\lg z_0'\rg$-block systems of $Q$, say
 ${\cal C_{w'}}$ on $W'$ and ${\cal C_{u'}}$ on $U'$, respectively.
Let $K_1$ and  $K_2$ be the respective  kernel of $Q$ acting on ${\cal C_{w'}}$ and ${\cal C_{u'}}$.
Since  $K_1\cup K_2\subsetneq Q$,  there exists an element $b\in  Q\setminus (K_1\cup K_2)$. Then $b$ moves $p$-blocks on both $W'$ and $U'$.
Therefore, $b^p$ fixes ${\cal C_{w'}}\cup  {\cal C_{u'}}$ pointwise and $b^{p^2}=1$. Set $Q_1=\lg z_0', b\rg $.  Then    $Q_1$
 is clearly transitive on both $W'$ and $U'$.  Since $Q_1$ is abelian and faithful on both sets, it acts regularly on both of them. Moreover,
 $Q_1\cong \ZZ_p^2$ if $b^p=1$; and  $Q_1\cong \ZZ_{p^2}$ if $b^p\ne 1$ (so $\lg b^p\rg =\lg z_0'\rg$).
  This completes the proof.
\qed

\begin{lemma} \label{main11} Using the above notation, suppose that $\op$ acts faithfully on both ${\cal B_w}$ and ${\cal B_u}$. Then $X$ admits a Hamilton cycle.
\end{lemma}
\demo
 Recall that $\op=P/K$ and $Z=Z(P)=\lg z_0\rg \leq K$.  Consider the action of $\op$ on ${\cal B_w}\cup {\cal B_u}$.
By Lemma~\ref{twofaithful}, there exists an abelian subgroup $H/K\leq \op$ that is regular on both ${\cal B_w}$ and ${\cal B_u}$. We will then complete the proof by considering two cases as follows.

\vskip 3mm
{\it Case 1: $\lg z_0\rg=K$:}
\vskip 3mm
Therefore, the subgroup   $Z(P)H$ has order $p^3$ and so it is regular on both $W$ and $U$. By Lemma~\ref{p3}, $X$ is Hamiltonian.
\vskip 3mm
{\it Case 2: $\lg z_0\rg\lneq K$:}
\vskip 3mm
In this case, we have $K_{(B)}\neq 1$, for any $B\in {\cal B_w}\cup {\cal B_u}.$
Since $H/K$ acts regularly on both $\mathcal{B}_w$ and $\mathcal{B}_u$ and $|\mathcal{B}_w| = |\mathcal{B}_u| = p^2$, it follows from Propositions~\ref{coset-graph-pro} and~\ref{2p^2DM} that there exists a Hamilton cycle in the quotient graph $X_{\langle z_0 \rangle}$, where $V(X_{\langle z_0 \rangle}) = \mathcal{B}_w \cup \mathcal{B}_u$.
 Then $X$ has a Hamilton cycle, by Lemma~\ref{main-lemma}.
\qed

\vskip 3mm
\subsection{$\op$ acts unfaithfully on both ${\cal B_w}$ and ${\cal B_u}$}
 Again we first present a group-theoretic lemma.

\begin{lemma}\label{twounfaithful}
Let $Q$ be a $p$-group that admits two unfaithful transitive representations of degree $p^2$ on $W'$ and $U'$, respectively.
 Suppose $Q$ is faithful on $W'\cup U'.$
   Then $Q$ contains a minimal subgroup $Q_1$ that is transitive on both $W'$ and $U'$,
where either $Q_1\cong\ZZ_p^2$ or $\ZZ_{p^2}$, acting regularly on both sets; or
 $Q_1=\lg a\rg \times \lg b\rg \cong \ZZ_{p}\times\ZZ_{p^2}$, acting unfaithfully  on both sets with the respective kernel $\lg a\rg $and $\lg b^p\rg$.
\end{lemma}
\demo
Let $K_{(W')}$ and $K_{(U')}$ be the kernels of $Q$ acting on $W'$ and $U'$, respectively.
By Proposition~\ref{pcentral}, we can take $z_1\in K_{(W')}\cap Z(Q)$ and $z_2\in K_{(U')}\cap Z(Q)$ such that $o(z_1)=o(z_2)=p$.
Then $\lg z_1\rg$  (resp. $\lg z_2\rg$ )  induces a $p$-block system of $Q$, say
 ${\cal C_{u'}}$ on $U'$ (resp. ${\cal C_{w'}}$ on $W'$). Let $K_1$ and $K_2$ be the respective kernels of $Q$ acting on these block systems. Clearly, the central subgroup  $\lg z_1z_2\rg$ is regular on every such $p$-block in ${\cal C_{w'}}\cup {\cal C_{u'}}$. Since  $K_1\cup K_2\subsetneq Q$,  there exists an element $b\in  Q\setminus (K_1\cup K_2)$. Then $b$ moves  $p$-blocks on both $W'$ and $U'$ and satisfies $b^{p^2}=1$.
 Thus we have the following two cases:
\vskip 3mm
(1) $o(b)=p$: in this case,  the abelian group $Q_1:=\lg z_1z_2, b\rg \cong \ZZ_p\times\ZZ_p$, acting regularly on both sets;
\vskip 3mm
(2) $o(b)=p^2$: since $b$ acts faithfully on $W'\cup U'$, the group  $\lg b\rg $ acts transitively at least one bipart. First suppose that  $\lg b\rg $ acts transitively on both $W'$ and $U'$. Then  $Q_1:=\lg b \rg$  is   regualr on both sets;
Secondly, $\lg b\rg $ is transitive on one bipart, say $W'$ but intransitive on $W'$. Then   set $a=z_1$ so that  the abelian group $Q_1:=\lg a\rg \times \lg b\rg\cong \ZZ_{p}\times \ZZ_{p^2}$, acts unfaithfully   on both sets with the respective kernel $\lg a\rg $ and $\lg b^p\rg $.
\qed

\begin{lemma} \label{main12} Using the above notation, suppose that $\op$ acts unfaithfully on both ${\cal B_w}$ and ${\cal B_u}$. Then $X$ admits a Hamilton cycle.
\end{lemma}
\demo
Recall that $\op=P/K$, where $K=K_w\cap K_u$ and $Z=Z(P)=\lg z_0\rg\leq K$.
Since $K$ is the kernel of $P$ acting on the set ${\cal B_w}\cup {\cal B_u}$, we have $x^p=1$ for any $x\in K^{V(X)}.$

 Suppose that $\op$ acts unfaithfully on both ${\cal B_w}$ and ${\cal B_u}$.
 By Lemma~\ref{twounfaithful}, there exists an abelian subgroup $\overline{Q_1}\leq\op$ that acts transitively on both ${\cal B_w}$ and ${\cal B_u}$; moreover, either $\overline{Q}_1$ acts faithfully and regularly on both ${\cal B_w}$ and ${\cal B_u}$, or $\overline{Q}_1\cong\ZZ_{p^2}\times\ZZ_p$ acts unfaithfully (while still transitive) on both ${\cal B_w}$ and ${\cal B_u}$. We now complete the proof in two steps.

\vskip 3mm
{\it Step 1: Find a Hamilton cycle for the quotient graph $X_{Z}$.}
\vskip 3mm
Suppose that $\overline{Q}_1$ is faithful and regular on both ${\cal B_w}$ and ${\cal B_u}$. Then the quotient graph $X_{Z}$ has a Hamilton cycle, by Propositions \ref{coset-graph-pro} and \ref{2p^2DM}.

Suppose that  $\overline{Q}_1=\lg \olb,\ola\rg\cong\ZZ_{p^2}\times\ZZ_p$ is unfaithful and transitive on both ${\cal B_w}$ and ${\cal B_u}$.
Followed from the proof in Lemma \ref{twounfaithful}, consider the action of $\overline{Q}_1$ on the quotient graph $(X_{Z})_{\lg\ola\olb^p\rg}$ of the graph $X_{Z}.$
Then $\lg\olb^p,\ola\rg\cong\ZZ_p\times\ZZ_p$ is in the kernel of $\overline{Q}_1$ acting on the set $V((X_{Z})_{\lg\ola\olb^p\rg})={\cal C_w}\cup{\cal C_u}$, where ${\cal C_w}$ (resp. ${\cal C_u}$) is the set of $\lg\ola\olb^p\rg$-orbits on ${\cal B_w}$ (resp. ${\cal B_u}$).  Observe that $\lg\olb\rg/\lg \olb^p\rg$ is regular on both ${\cal C_w}$ and ${\cal C_u}$.
Then, by Propositions \ref{coset-graph-pro} and \ref{2p^2DM}, $(X_{Z})_{\lg\ola\olb^p\rg}$ has a Hamilton cycle. Note that $\lg\olb^p,\ola\rg\cong\ZZ_p\times\ZZ_p$ lies in the kernel of $\overline{Q}_1$ acting on the set $V((X_{Z})_{\lg\ola\olb^p\rg})={\cal C_w}\cup{\cal C_u}$.  Since each block in ${\cal C_w}\cup {\cal C_u}$ has size $p$, Lemma \ref{main-lemma} implies that $X_{Z}$ contains a Hamiltonian cycle.

\vskip 3mm
{\it Step 2:  Find a Hamilton cycle for $X$.}

\vskip 3mm
We consider two cases, separately.
\vskip 3mm
{\it Case (i)  $K=Z(P):$}
\vskip 3mm
Suppose that $\overline{Q}_1=\lg \olb,\ola\rg$ is faithful and regular on both ${\cal B_w}$ and ${\cal B_u}$. Then $Z(P)Q_1$ is regular on $W$ and $U$. By Lemma \ref{p3},  $X$ has a Hamilton cycle,

Suppose that  $\overline{Q}_1=\lg \olb,\ola\rg\cong\ZZ_{p^2}\times\ZZ_p$ is unfaithful and transitive on both ${\cal B_w}$ and ${\cal B_u}$. Note that $\overline{Q}_1=\lg \olb,\ola\rg$, where $[a,b]=z_0^i$. This implies that $[a,b^p]=1$, and so $b^p\in Z(Q_1),$ where $Q_1/K=\overline{Q_1}$.
 Without loss of generality, assume that $\olb^p$ fixes every $\lg z_0\rg$-orbit on $U$.
Now consider the action of $\lg b^p,z_0\rg\leq Z(\lg b,z_0, a\rg)$ on each $\lg z_0\rg$-orbit of $U$. Since $\langle a, b, z_0 \rangle$ acts transitively on $U$, there exists $u \in U$ such that the stabilizer $\langle b^p, z_0 \rangle_u$ is nontrivial and lies in the kernel of the action of $Q_1$ on $U$. This contradicts to the faithfulness of the action of $\langle b, z_0, a \rangle$ on $U$.

\vskip 3mm
{\it Case  (ii)  $Z\lneq K:$}
\vskip 3mm
In this case, we have $K_{(B)}\neq 1$ for any $B\in{\cal B_w}\cup {\cal B_u}.$ By Lemma \ref{main-lemma}, a Hamilton cycle of $X_Z$  can be lifted  a Hamilton cycle of $X$. \qed

\subsection{$\op$ acts unfaithfully on exactly one of ${\cal B_w}$ and ${\cal B_u}$}
Suppose that $\op $ acts unfaithfully on one bipart, say  ${\cal B_w}$ with the kernel $Z_1$ and faithfully on ${\cal B_u}$.
 Then $\op \le G=\ZZ_p\wr \ZZ_p$ and $|\op|\ge p^3$.    View $G$ with a affine group, that is $G=V\rtimes \lg \rho \rg $, where $V$ is identified with $V(p,p)$, the translation group, where
$v^\rho=(x_0, x_1, \cdots, x_{p-1})$, if $v=(x_1, x_2, \cdots, x_0).$ For convenience, we denote by $G_k := [G, G, \dots, G]$ the $k$-th term of the lower central series of $G$.
 To begin, we present a lemma of a group-theoretic nature. It is worth noting that $\ZZ_p \wr \ZZ_p$ is a Sylow $p$-subgroup of $S_{p^2}$.
\begin{lemma}\label{oneunfaithful}
Let $Q$ be a $p$-group  having a unfaithful transitive representations of degree $p^2$ on $W'$ and  a faithful transitive representations of degree $p^2$  on  $U'$.
   Then $Q$ contains a transitive subgroup  $Q_1=G_{k}\lg v\rho \rg $ acting faithfully  on  $U'$   and unfaithfully  on  $W'$  with the kernel $G_{k+1}$, where  $k\le p-1$ and $v\in V$.
\end{lemma}
\demo  Let $\rho(A)$ be the permutation matrix  corresponding to $\rho$ and $I$ the identify matrix.  Set $v_1=(1, 0\cdots, 0)$ and
$$
\begin{array}{lll}
v_i&:=&[v_1, \overbrace{\rho, \rho, \cdots , \rho}^{i-1}]= v_1(-I+\rho(A))^{i-1}=v_1(\sum\limits_{k=0}^{i-1}(-1)^{i-1-k}\left(
          \begin{array}{c}
            i-1  \\
            k
          \end{array}
        \right)\rho(A)^k)\\
&=&((-1)^{i-1},
 (-1)^{i-2}\left(
          \begin{array}{c}
            i-1  \\
            1
          \end{array}
        \right),
 (-1)^{i-3}\left(
          \begin{array}{c}
            i-1  \\
            2
          \end{array}
        \right),\cdots,
        \left(
          \begin{array}{c}
            i-1  \\
            i-1
          \end{array}
        \right),0,0,\cdots),\\
\end{array}$$
with $1\leq i\leq p.$
Note that $G=\lg v_1,\rho\rg$.
Then the derived group
$$
\begin{array}{lll}
G_2&=&\lg v_2, G_3\rg =\lg v_2^{\lg \rho\rg }\rg =\lg (x_{p-1}-x_0, x_0-x_1, \cdots , x_{p-2}-x_{p-1})\di  x_i\in \ZZ_p \rg\\
& =&\{v\in V\di \sum x_i=0\},
\end{array}$$
and more generally, the $i$-th lower central subgroup
$G_i=\lg v_i, G_{i+1}\rg =\lg v_i^{\lg \rho \rg}\rg $, where $2\le i\le p$.
In particular, $Z(G)=G_p=\lg  (1, 1,  \cdots , 1)\rg $, as $(-1)^\imath\left(
          \begin{array}{c}
            p-1  \\
            \imath
          \end{array}
        \right)=1$ for any $0\leq \imath\leq p-1,$
        and so the nilpotent class   $c(G)$ of $G$ is $p$.
Moreover,  $(v\rho)^p=1$ or  $(v\rho)^p\in Z(G)\setminus \{1\}$    if and only if   $v\in G_2$ or $v\in V\setminus G_2$, respectively. Finally, $G_i$ are all normal subgroups of $G$ of index more than $p$.
 Remind that some of the above information  are also available in \cite[Kapitel III]{Hup}.

Let $Q$ be a $p$-group  having a unfaithful transitive representations of degree $p^2$ on $W'$ and a faithful transitive representations of degree $p^2$ on $U'$.
Let $\lg z_0'\rg\leq Z(Q)$, where $o(z_0')=p$.
Consider the action of $Q$ on the set ${\cal C_u}$ of $\lg z_0'\rg$-orbits on $U'$. Then $Q/Q_{({\cal C_u})}\cong \ZZ_p$, and we may write $Q=(V\cap Q)\lg v\rho \rg $.
Furthermore, since $V \cap Q$ is normal in $G$, we have $V \cap Q =G_{l}$, for some integer $l$ so that $Q=G_{l} \langle v \rho \rangle$.

  Suppose that $k$ is the smallest  one such that $G_{k+1}$ fixes $W'$ pointwise but  $G_k$ does not.
Set  $Q_1=G_k \lg v\rho \rg $. Then  $Q_1$  acts faithfully on $U'$ and unfaithfully on $W'$  with the kernel $G_{k+1}$, where  $k\le p-1$ and $v\in V$.
\qed

\begin{lemma} \label{main13} Using the above notation, suppose that $\op $ acts unfaithfully on one bipart, say  ${\cal B_w}$  and faithfully on ${\cal B_u}$.   Then $X$ admits a Hamilton cycle.
\end{lemma}
\demo
Recall that $\op=P/K$, where $K=K_u\lneq K_w$ and $Z=Z(P)=\lg z_0\rg\leq K$.
Take an edge $wu$ and two $Z$-blocks ${\bf w}_0$ and ${\bf u}_0$ containing $w$ and  $u$, respectively.
Suppose $\op$ acts unfaithfully on one bipart ${\cal B_w}$ and faithfully on ${\cal B_u}$.
 By Lemma~\ref{oneunfaithful}, $\op$ contains
a subgroup $\op_1=G_k\lg \olb\rg $  of $\op$, where $\olb=v\rho$ for some $v\in V(p,p)$.
Then  either $\olb^p=1$ or $\olb^p=(x, x, \cdots, x)$ for some $x\in\ZZ_p$.
 Furthermore, $P_1=K.\op_1$ is  an extension of $\op_1$ by $K$, which
   acts transitively on both $W$ and $U$.
   We now consider two cases, separately.

  \vskip 3mm
 (1) $K=\lg z_0\rg$.
 \vskip 3mm
The arguments is divided into the following two steps.

\vskip 3mm
{\it Step 1: Analysis of the group structure.}
\vskip 3mm
Consider the transitive action of $P_2:=K.G_k$ (a group extension) on the orbits $w^{P_2}$ and  $u^{P_2}$. Each orbit has length $p^2$, and $|(P_2)_w|=|(P_2)_u|=p^{p-k}$, since $|G_k|=p^{p-k+1}$.
 Set $P_0=(P_2)_w\cap (P_2)_u$.
Consider the action of $P_2$ on $w^{P_2}$ and $u^{P_2}$, and the action of $G_k$ on the sets ${\cal B_w}$ and ${\cal B_u}$.
 Recall that ${\bf w}_0$ (resp. ${\bf u}_0$) is the $Z$-orbit containing $w$ (resp. $u$).
 Note that both ${\cal B_w}$ and ${\cal B_u}$ are $p$-block systems of $P$ on the set $W$ and $U$, respectively.
Then $(P_2)_w\cong((P_2)_wK)/K\leq (G_k)_{{\bf w}_0}$ and $(G_k)_{{\bf w}_0}=G_{k+1}$, as $G_{k+1}\leq (G_k)_{{\bf w}_0}$ and $|(G_k)_{{\bf w}_0}|=|G_{k+1}|$.
Hence, $(P_2)_w\cong\overline{(P_2)_w}=(\overline{P_2})_{\bf w_0}= G_{k+1}$, as $|G_{k+1}|=|(P_2)_w|=p^{p-k}.$
Note that $G_{k+1}$ is normal in $\op_1$ and $\op$ is faithful on ${\cal B_u}$, we get that $G_{k+1}$ can not fix any vertex in ${\cal B_u}$. 
This implies that
 $|(P_2)_w:P_0|=p$, which yields  $|P_0|=p^{p-k-1}$. Moreover,  we have
 $$(P_2)_{w}\cong (\op_2)_{{\bf w_0}}\cong \ZZ_p^{p-k} \,  \,{\rm and}\, \,  (P_2)_{u}\cong (\op_2)_{{\bf u_0}}\cong \ZZ_p^{p-k}.$$
Therefore,  there exit $a, z\in P_2$ such that
$(P_2)_{w}=P_0\times\lg z\rg $ and $ (P_2)_{ u}= P_0\times\lg a\rg $, which implies
$$P_2=(P_0\times\lg z_0\rg) \lg a, z\rg ,$$
where $[P_0\times\lg z_0\rg, \lg a, z\rg ]=1$ and either $[a, z]=z_0$  or   $[a, z]=1$.
Since $(\op_2)_{{\bf w_0}}=\overline{P_0}\times\lg \olz\rg $ and $(\op_2)_{{\bf u_0}}=\overline{P_0}\times\lg \ola\rg $,
one may choose $\olz=(1, 1, \cdots , 1)\in G_{k+1}$ and $\ola=v_k^{\rho}=(0, \ast, \cdots, \ast)\in G_k\setminus G_{k+1}$ so that $\lg z, a, z_0\rg $ is transitive on both  $w^{P_2}$ and  $u^{P_2}$.

Let $b$ be a perimage of $\olb$ and set
 $$P_1=P_2 \lg b\rg .$$
 Suppose that $[a, z]=z_0$.  Then $Z(P_2)=P_0\times \lg z_0\rg$ is a characteristic subgroup of $P_2$.
  Since $b$ normalizes $Z(P_2)$, we have that $\olb$ normalizes $\op _0\cong P_0$. Therefore, $\op_0\lhd G$, where $G\cong\ZZ_p\wr\ZZ_p$.
 Since $G_i$ ($2\le i\le p$) are all normal subgroups of $G$ with order $p^{p+1-i}$,
   we get  $\lg \olz \rg \le \op_0$, a contradiction, as $\op_0 \lg \olz \rg\cong P_0\times\lg z\rg $. Therefore, we have $[a,z]=1$.

Let $C=C_{P_1}(\lg z_0, z\rg)$. Then $C\ne P_1$, otherwise, $z\in Z(P_1)$ and so $z$ fixing  $W$ pointwise, as $z\in (P_1)_w$.
Since  $P_1/C\le \Aut(\lg z_0, z\rg)\cong \GL(2,p)$, we have $|P_1:C|=p$. Then $C=P_2$, so that $[b, z]\ne 1$, as $[a,z]=1$.
 Since $\olb^p=1$ or $\olb^p\in\lg\olz\rg$, we write
 $b^p=z^iz_0^j$ for some $i,j\in\ZZ_p$.
 If $z^i\ne 1$, then   $[z^i, b]=[b^pz_0^{-j}, b]=1$ and so $[z, b]=1$, a contradiction. Therefore, we have that
 $b^p=z_0^j$ for some $j\in\ZZ_p$.

\vskip 3mm
{\it Step 2: Find a Hamilton cycle of $X$.}
\vskip 3mm
Recall that $P_1=(\lg z_0\rg \times \lg P_0, a, z\rg \rg )\lg b\rg $, where $\lg P_0, a, z\rg \cong G_k$, $[a,z]=1$ and $o(a)=o(z)=p$.
Set ${\bf w}=w^{\lg az\rg}$ and  ${\bf u}=u^{\lg az\rg}$.
 In  viewing of
$$\begin{array}{l}(P_1)_{w}=P_0\times\lg z\rg\lneqq (P_1)_{w}\times \lg az\rg  \lneqq((P_1)_{w}\times \lg az\rg )\times \lg z_0\rg)\lg b\rg ,\\
(P_1)_{u}=P_0\times\lg a\rg\lneqq (P_1)_{u}\times \lg az\rg \lneqq  ((P_1)_{u}\times \lg az\rg )\times \lg z_0\rg)\lg b\rg ,\end{array}$$
 we get two complete block systems on $W$ and $U$:
 $${\cal W}=\{ {\bf w}^{z_0^ib^j}\mid i,j\in \ZZ_p\} \quad {\rm and}\quad {\cal U}=\{ {\bf u}^{z_0^ib^j}\mid i,j\in \ZZ_p\}.$$
\f Let $Y_1$ be the corresponding block graph.
We now prove that for any two adjacent blocks $B_1,B_2\in V(Y_1)$, the induced graph $X[B_1,B_2]\cong K_{p,p}.$
Without loss of generality, assume that $w'$ is adjacent to $u'$, where $w'\in B_1$ and $u'\in B_2$. Then $(P_2)_{w'}K/K=G_{k+1}=(P_2)_{w}K/K$, it follows from $z\in (P_2)_{w}$  that $z\in (P_2)_{w'}$. Because $z$ acts transitively on $B_2$, we obtain $X[B_1, B_2] \cong K_{p,p}$, as desired.

Note that $b^p=z_0^j$ for some $j\in\ZZ_p$. Then $\lg z_0, b\rg \cong \ZZ_p^2$ or $\ZZ_{p^2}$, and so $\lg z_0, b\rg$ acts regularly on both ${\cal W}$ and  ${\cal U}$.
By Propositions \ref{coset-graph-pro} and \ref{2p^2DM},
$Y_1$ has a Hamilton cycle, say $C_1$.
Since $X[B_1, B_2]\cong K_{p,p}$ for any two adjacent blocks $B_1\in{\cal W}$ and $B_2\in{\cal U}$,  $C_1$ can be lift to a Hamilton cycle of $X$ by Proposition \ref{pro:4}.

 \vskip 3mm
   (2) $\lg z_0\rg \lneqq K$.
    \vskip 3mm
    Consider the  quotient graph $Y=X_K$. Then $\op_1=G_k\lg \olb\rg $. In our case,  $G_k$ induces $p$-blocks on both ${\cal B_w}$ and ${\cal B_u}$, say ${\cal W}$ and ${\cal U}$,  and by $Y_{G_k}$ we denote the corresponding quotient graph of $Y$. Since $Y_{G_k}$ is of order $2p$, with a subgroup of order $p$ acting transitively on both biparts,
    we may take a cycle $\mathcal{C}_2$ by Propositions \ref{coset-graph-pro} and \ref{2p^2DM}.
    Since $\op_{{\bf w}_0}$ and $\op_{{\bf u}_0}$ are different, any induced subgraph by two adjacent $p$-blocks in   ${\cal W}$ and ${\cal U}$ respectively  is isomorphic to $K_{p,p}$. Clearly, $\mathcal{C}_2$ can be lifted a Hamilton cycle $\mathcal{C}_1$ of $Y$ by Proposition \ref{pro:4}. Since  $\lg z_0\rg \lneqq K$, we have $K_{(B)}\neq 1$ with $B=w^{\lg z_0\rg}$. By Lemma \ref{main-lemma},  $\mathcal{C}_1$ can be lifted a Hamilton cycle $\mathcal{C}$ of $X$.\qed

\section{$A$ acts unfaithfully on one bipart}
 The main result of this section is the following theorem.
\begin{theorem}\label{main2}
 Using the above notation,
suppose that  $A$ acts unfaithfully on one bipart. Then $X$ has a Hamilton cycle.
\end{theorem}

Suppose $A$ acts  unfaithfully on one bipart, say $W$.
Then all graphs $X$ have already been classified by Du and Wang in \cite{Du-Wang}. Based on this classification, Theorem \ref{main2} is an immediate consequence of Lemmas \ref{12}, \ref{34} and \ref{58}.
\vskip 3mm
To state their classification, we recall  the following  definition.

\begin{definition} Let $\Sigma=(V_1,E_1)$ be a connected edge-transitive graph with bipartition $V_1=W_1\cup U_1$, where $|W_1|=p^3$ and $|U_1|=p^2$ for an odd prime $p$. Define a bipartite graph $\Gamma=(V,E)$ with bipartition $V=W\cup U$, where
$$
\begin{array}{l}
  W=W_1, U=\{(u,i)\mid u\in U_1, i\in\ZZ_p\}\,\,\text{and}\,\,
  E=\{\{w,(u,i)\}\mid \{w,u\}\in E_1, i\in\ZZ_p\}.
\end{array}
$$
Then $\Gamma$ is called the graph {\it expanded} from $\Sigma$.
\end{definition}
\f In fact, each automorphism $\phi_1$ of $\si$ can be lifted to an automorphism $\phi$ of $X$ as follows:
\begin{eqnarray}\label{Eq1} \phi (w)=\phi_1(w), \, \phi ((u, k))=(\phi_1(u), k),\end{eqnarray}
for any $w\in W=W_1$, $u\in U_1$ and $k\in \ZZ_p$. Moreover, choose a permutation $\s_1 $ on $Z_p$. Define an automorphism   $\s$ on $V(X)$ by
\begin{eqnarray}\label{Eq2}\s (w)=w\quad {\rm and}\quad \s(u,k)=(u, \s_1(k)),\end{eqnarray}
for any $w\in W$, $u\in U_1$, $k\in Z_p$.  Clearly, $[\s, \phi]=1$ for any $\phi$ lifted from $\phi_1\in \Aut(\si).$

Now we are ready to state the classification theorem in  \cite{Du-Wang}.
\begin{proposition}\label{unf-c}
 Let $p$ be an  odd prime and  let  $X$ be  a semisymmetric graph of order $2p^3$ whose automorphism group acts unfaithfully on one bipartite part. Then $X$ is expanded from one of the families of graphs $\si =\si_{ij}$ with $1\le i\le 8$, $i\ne 4$, $j=0, 1$; or from $\si_4$,     as will be given in the following four subsections.
 \end{proposition}

  \subsection{{\bf Families  1 and 2}}
First we present a bi-coset graph description of the graphs $\Sigma$ and then   prove  that all expanded graphs $X$ of such $\Sigma$ (from these families) admit a Hamilton cycle.
\vskip 3mm
 \f {\it {\bf Family 1:} Let   $M=\lg a, b\rg \cong S_3$, where  $o(a)=3$ and $o(b)=2$. Set
$$A=M\times M\times M: S_3,\,\,L=\lg b\rg\times\lg b\rg\times\lg b\rg:S_3,\,\,R= (\lg b\rg \times M\times M)\rtimes \lg (23)\rg,
$$
$D_1=RL$ and $D_2=RaL.$
 Define two bi-coset
graphs $\si _{11}={\bf B}(A; L, R, D_1)\, \, {\rm \,and }\,\, \si _{12}={\bf B}(A; L, R, D_2).$  The
valency $d(L)$ of the vertex $L\in [G:L]$ is $3$  for $\si_{11}$ and $6$ for $\si_{12}$.
}

\vskip 3mm

 \f {\it {\bf Family 2:}} {\it Let $p\ge 5$ be an odd prime and let $\Lambda $ be the subgroup of $\FF_p^*$ with order $r\ge 1$.
  Set
 $$\begin{array}{lll}
 x&=&{\footnotesize \left(\begin{array}{ccc} 1 &2 &2 \\ 0&1 &2 \\
0&0 &1\end{array}\right)},
 \,\,  y=y(\lambda, \mu )={\footnotesize \left(\begin{array}{ccc}  \mu ^2\lambda^{-1} &0 &0 \\ 0&\mu &0\\
0&0 &\lambda\end{array}\right)} \in \GL(3,p),\\ \\ S&=&\lg y(\mu, \lambda )\di  \mu \in \FF_p^*, \lambda \in \Lambda \rg \le
N_{\GL(3,p)}(\lg x\rg).\end{array}$$
Let  $N=\lg t_1,t_2,t_3\rg$ and $N_0=\lg t_1, t_2\rg \le
N,$ where $t_1=(1, 0, 0), t_2=(0, 1, 0)$ and $t_3=(0, 0, 1).$
Set $A=N\rtimes (\lg x\rg \rtimes S), \,\, L=\lg x\rg \rtimes  S, \,\,
 R=N_0\rtimes S,$
 where $t_1^x=t_1t_2^2t_3^2$, $t_2^x=t_2t_3^2$ and $t_3^x=t_3.$
Define the graphs
$$\begin{array}{lll} \Sigma _{21}(p)&=&{\bf B}(A; L, R, D), \, {\rm \,where\,}  \Lambda=\FF_p^*, D=RL;\\
\Sigma _{22}(p, r)&=&{\bf B}(A; L, R, D), \, {\rm \, where\,}  2\le r,\, D=Rt_{3}L,\end{array}$$
where  $d(L)=p$ for $\Sigma_{21}(p)$ and $d(L)=rp$ for $\Sigma_{22}(p,r)$.  }

  \begin{lemma}\label{12}
  Suppose that  $\si\in \{\si_{11}$, $\si_{12}$,
$\Sigma_{21}(p)$,  $\Sigma_{22}(p,r)$\}. Then $X$ admits a Hamilton cycle.
\end{lemma}
  \demo (1)  $\si=\Sigma_{21}(p)$ or $\Sigma_{22}(p,r)$, where $p\ge 5$
    \vskip 3mm
    Correspondingly, we set $X = X_{21}$ and $X = X_{22}$ for $\Sigma=\Sigma_{21}(p)$ or $\Sigma_{22}(p,r)$, respectively.
     Let $P_1=N\rtimes \lg x\rg $ be a Sylow $p$-subgroup of $A$.  Set $w=L$ and $u=R$ so that  two biparits of $\si$  can be written by
  $W_1=\{ w^g\di g\in P_1\} \quad {\rm and} \quad U_1=\{ u^g\di g\in P_1\}.$ Let $P$ be the lifts of $P_1$ as defined  in Eq(\ref{Eq1}), where we identify $P$ with $P_1.$
  Then  two biparits of $X$ are
   $$W=W_1=\{ w^g\di g\in P\}  \quad {\rm and} \quad U=\{ (u, k)^g\di g\in P, k\in \ZZ_p\}.$$
  Let $\s$ be  as defined  in Eq(\ref{Eq2}) so that $[\s, P]=1$.

Noting  $t_1=(1, 0, 0), t_2=(0, 1, 0)$, $t_3=(0, 0, 1)$, where $t_3\in Z(P)$ so that $\lg t_3, xt_2\rg \cong \ZZ_p^2$.  Let $$\begin{array}{lll} &{\bf w}_{0}=\{w^g\di g\in \lg t_3, xt_2\rg \},\quad &{\cal W}=\{ {\bf w}_{0}^{t_1^i}\di i\in \ZZ_p\},\\
&{\bf u}_{0}=\{(u, 0)^g\di g\in \lg t_3,  xt_2\rg \},\quad &{\cal U}=\{ {\bf u}_{0}^{\s ^i}\di i\in \ZZ_p\},\end{array}$$
By the definition, $w\sim (u, k)$ for $X_{21}$ and  $w\sim (u^{t_3}, k)$ for $X_{22}$.  Let $X_0$ be the induced subgraph $X[{\bf w}_{0}, {\bf u}_{0}]$ of $X$. Consider the action  $P_0=\lg t_3, t_2, x\rg $ on  $X_0$. Then $(P_0)_{w}=\lg x\rg $ and $(P_0)_{(u, 0)}=\lg t_2\rg $. Since $\lg DD^{-1}\rg =\lg (P_0)_{w}, (P_0)_{(u, 0)}\rg =P_0$ for both $X_{21}$ and  $X_{22}$, by Proposition \ref{coset-graph-pro},
$X_0$ is a connected bi-Cayley graph of $ \lg t_3, xt_2\rg \cong \ZZ_p^2$.  By Propositions \ref{coset-graph-pro} and \ref{2p^2DM} , it  has a Hamilton cycle, say ${\cal C_0}=(u, 0), \cdots , w, (u, 0)$ for $X_{21}$, and ${\cal C_0}=(u^{t_3}, 0), \cdots , w, (u^{t_3}, 0)$ for $X_{22}$.
   Let ${\cal P_0}=(u, 0), \cdots , w$ (resp. ${\cal P_0}=(u^{t_3}, 0), \cdots , w$) be the corresponding  Hamilton path of ${\cal C_0}$ by removing  $(u, 0)$ (resp. $(u^{t_3},0)$) for $X_{21}$ (resp. $X_{22}$). Since $\lg t_1\s\rg $ is respectively transitive on  blocks  of ${\cal W}$  and  blocks of ${\cal U}$,   $w\sim (u, 0)^{t_1\s}=(u, 1)$ for $X_{21}$  and  $w\sim (u^{t_3}, 0)^{t_1\s}=(u^{t_3}, 1)$ for $X_{22},$
    we get a $H$-cycle of $X_{21}$:
    $$\mathcal{P}_0\sim\mathcal{P}_1\sim\cdots\sim\mathcal{P}_{p-1}\sim (u,0);$$
 and that of $X_{22}:$
  $$\mathcal{P}_0\sim\mathcal{P}_1\sim\cdots\sim\mathcal{P}_{p-1}\sim (u^{t_3},0),$$
  where
  $\mathcal{P}_i=(u, 0)^{(t_1\sigma)^i}, \cdots , w^{(t_1\sigma)^i}$ for $X_{21}$
  and that
  $\mathcal{P}_i=(u^{t_3}, 0)^{(t_1\sigma)^i}, \cdots , w^{(t_1\sigma)^i}$ for $X_{22}$.
\vskip 3mm
(2)  $\si=\Sigma_{11}$ or $\Sigma_{12}$, where $p=3$.
    \vskip 3mm
Correspondingly, set $X=X_{11}$ and $X_{12}$ for $\Sigma=\Sigma_{11}$ or $\Sigma_{12}$, respectively. Let $P_1=N\rtimes \lg x\rg \le \AGL(3,3)$ be a Sylow $p$-subgroup of $A$, where $N=V(3,3)$ and $t^x=(x_0, x_1, x_2)$ for ant $t=(x_1, x_2, x_0)\in N$.
        No loss, let $(P_1)_w=\lg x\rg $ and  $(P_1)_u=\lg (0, x_1, x_2)\di x_1, x_2\in \ZZ_3\rg $. Then we have the completely same arguments as used in (1),  by resetting
        $t_1=(0, 2, 2),\, t_2=(0, 1, -1)$ and   $t_3=(1,1,1)$ so that $t_3\in Z(P)$,
        $t_1^x=t_1t_2^2t_3^2$, $t_2^x=t_2t_3^2$ and $t_3^x=t_3.$
 \qed

  \subsection{Families 3 and 4}
\f {\it {\bf Family 3:}}  {\it For $p\ge 5$, let $P=\lg a, b, c\di a^p=b^p=c^p=1, [b, a]=c, [c,a]=[c,b]=1\rg .$
 For any $s', t'\in \FF_p^*=\lg \th \rg $, set $\phi(s', t')\in \Aut(P)$ such that $a\mapsto a^{s'}, b\mapsto b^{t'}$.
 Let     $r_1$ and $r_2$ be  two divisors  of $p-1$ with $r_2\ge 3$. Define
 an abelian subgroup $S=\lg \phi(1, 2), \phi(s, t)\rg $ of order $r_1r_2$ such that
 $S\cap \lg \phi(1, \th)\rg =\lg \phi(1, e)\rg \cong Z_{r_1}, \, S/\lg \phi(1, e)\rg\cong Z_{r_2},$ that is $o(e)=r_1, o(s)=r_2$ and $t^{r_2}\in \lg e\rg$. Now, $T=\lg t\in \FF_p^*\di t^{r_2}\le \lg e\rg \rg $ is a subgroup of $\FF_p^*$  of order $r_1(r_2, \frac{p-1}{r_1})$ and so  we choose
  $ e=\th^{\frac{p-1}{r_1}},\,   s=\th^{\frac{p-1}{r_2}},\,  t=(\th^{\frac{p-1}{r_1(r_2,
\frac{p-1}{r_1})}})^{r_3},$
 where $0\le r_3\le  (r_2, \frac{p-1}{r_1})-1$, and if $r_1=1$ then we let $r_3\not\in \{  1, r_2-1\}$.
Let
 $$A=P\rtimes S, \quad L=S, \quad R=\lg b\rg S.$$
Define the graphs $$\begin{array}{lll} \si _{31}(p, r_2)&=&{\bf B}(A; L, R, D),   D=RaL, \, {\rm \,where\,}  r_1=p-1,
r_3=0;\\ \si _{32}(p, r_1, r_2, r_3)&=&{\bf B}(A; L, R, D), \,  D=RcaL,\end{array}$$
where  $d(L)=r_2$ for $\si_{31}(p,r_2)$ and $d(L)=r_1r_2$ for $\si _{32}(p, r_1, r_2, r_3)$.
  }
\vskip 3mm
 \f {\it {\bf Family 4:}} {\it Let $p\ge 5$ be an odd prime and let $P=\lg a, b\di a^{p^2}=b^p=1, [b, a]=a^{p}\rg .$
Let $\FF_{p^2}^*=\lg \theta\rg$ and  $r\ge 3$,  a divisor  of $p-1$. Let
  $S=\lg \phi(\theta^{\frac{p(p-1)}r})\rg \le \Aut(P)$, where $\phi (a)=a^{\theta^{\frac{p(p-1)}r}}$ and  $\phi(b)=b$. Let $$A=P\rtimes S, \quad L=S,\quad R=\lg b\rg S,\quad D=RaL.$$ Define a  graph
$\si_4(p, r)={\bf B}(A; L, R, D)$, where  $d(L)=r$.
  }

\begin{lemma}\label{34}
 Suppose that  $\si=\Sigma_{31}(p,r_2)$, $\Sigma_{32}(p,r_1,r_2,r_3)$ or $\si_4(p, r)$. Then $X$ admits a Hamilton cycle.
\end{lemma}
  \demo (1) Let $\si=\Sigma_{31}(p,r_2)$ or $\Sigma_{32}(p,r,r_2,r_3)$. Correspondingly, set $X = X_{31}$ and $X = X_{32}$, respectively. Recall that $P_1:=\lg a, b, c\di a^p=b^p=c^p=1, [b, a]=c, [c,a]=[c,b]=1\rg  $ is a Sylow $p$-subgroup of $A$.  Set $w=L$ and $u=R$ so that  two biparits of $\si$  can be written by
  $W_1=\{ w^g\di g\in P_1\} \quad {\rm and} \quad U_1=\{ u^g\di g\in P_1\}.$ Let $P$ be the lifts of $P_1$ as defined  in Eq(\ref{Eq1}), where we identify $P$ with $P_1.$
  Then  two biparits of $X$ are
   $$W=W_1=\{ w^g\di g\in P\}  \quad {\rm and} \quad U=\{ (u, k)^g\di g\in P, k\in \ZZ_p\}.$$
  Let $\s$ be  as defined  in Eq(\ref{Eq2}) so that $[\s, P]=1$.
Then the subgroup $P':=\lg a,c,b\sigma\rg\cong P$ is regular on both biparts $W$ and $U$. It followed from Lemma \ref{p3} that the graph $X$ admits a Hamilton cycle.

\vskip 3mm (2)  Let $\si=\Sigma_{4}(p,r)$,  and set  $X=X_{4}$. Recall that $P_1=\lg a, b\di a^{p^2}=b^p=1, [b, a]=a^{p}\rg $ is  a Sylow $p$-subgroup of $A$. Then we have the completely
same arguments as that in (1), just noting that the subgroup $P':=\lg a,b\sigma\rg\cong P$ is regular on both biparts $W$ and $U$.
\qed

\subsection{Families 5-8}

Following \cite{WDL}, families 5-8, that is, $\Sigma_{5i}$, $\Sigma_{6i}$, $\Sigma_{7i}(p)$, and $\Sigma_{8i}(p,r)$ with $i=1,2$, give rise to the following graphs.
\begin{definition}\label{X} Define six families of bipartite graphs $X$ with bipartition $V(X)=W(X)\cup U(X)$, where $W(X)=\{(i,j,k)\mid i,j,k\in\ZZ_p\},\,\,\,\,
U(X)=\{[x,y,z]\mid x,y,z\in\ZZ_p\},$
and edge set
$$
\begin{array}{lll}
E(X)&=&\{\{(i,j,k),[x,i+b,k+\frac{p-1}{2}]\}\mid i,j,k,x\in\ZZ_p,b\in \Lambda\}\cup\\
&&\{\{(i,j,k),[x,j+sb,k+\frac{p+1}{2}]\}\mid i,j,k,x\in\ZZ_p,b\in \Lambda\},
\end{array}
$$
where $s=\theta^{\frac{p-1}{2r}}$, $\ZZ_p^*=\lg\theta\rg$ for the following subfamily of graphs $X_2(p,r),$ $s=1$ for other five subfamilies of graphs $X_i(p,r)$, and $\Lambda$ is given by
\begin{enumerate}
  \item {\it Graphs $X_1(p,r):$} Let $p\geq3$ and let $\Lambda$ be a subgroup of $\ZZ_p^*$ of order $r$, where $(p,r)\neq (7,3),(11,5).$
  \item {\it Graphs $X_2(p,r):$} Let $p\geq5$ and let $\Lambda$ be a subgroup of $\ZZ_p^*$ of order $r\geq 2$, where $(p,r)\neq (7,3),(11,5)$ and $2r\mid (p-1)$.
  \item {\it Graphs $X_3(11,5):$} Let $p=11$ and $\Lambda=\{0,2,3,4,8\}\subseteq\ZZ_{11}$.
  \item {\it Graphs $X_4(11,6):$} Let $p=11$ and $\Lambda=\{1,5,6,7,9,10\}\subseteq\ZZ_{11}$.
  \item {\it Graphs $X_5(p,r):$} Choose a point $\lg v\rg$ and a hyperplane $\mathcal{L}$ in the project space $\PG(n-1,q)$, where $\frac{q^n-1}{q-1}=p\geq 7$, and let $G=\lg t\rg$ be a Singer subgroup of $\PGL(n,q)$. Let $\Lambda=\{l\in\ZZ_p\mid \lg v\rg\in\mathcal{L}^{t^l}\},$ where $r=|\Lambda|=\frac{q^{n-1}-1}{q-1}.$
  \item {\it Graphs $X_6(p,r):$} Adopting the same notation as in (5), set $\Lambda=\{l\in\ZZ_p\mid \lg v\rg\notin \mathcal{L}^{t^l}\},$ where $r=q^{n-1}$.
\end{enumerate}
\end{definition}

Using the definitions of the graphs $\Sigma$ from Families 5-8, we show in Lemma \ref{58} that all expanded graphs $X$ of such $\Sigma$ admit a Hamilton cycle.

\begin{lemma}\label{58}
Let $X = X_i(p,r)$ be as defined in Definition \ref{X}, where $1 \leq i \leq 6$. Then $X$ admits a Hamilton cycle.
\end{lemma}
\demo
Now we divide our proof into two cases.
\vskip 3mm
{\it Case 1: $s=1$.}
\vskip 3mm
With the above notations, we define a quotient graph $X_{\lg \rho\rg}$ which is induced by the automorphism $\rho:(i,j,k)\mapsto (i,j,k+1)$
and $[x,y,z]\mapsto [x,y,z+1]$, where $i,j,k,x,y,z\in\ZZ_p$ and $\lg \rho\rg\cong\ZZ_p$.
Then we get a quotient graph $X_{\lg \rho\rg}$ induced by $\lg \rho\rg$, where $V(X_{\lg \rho\rg})=W(X_{\lg \rho\rg})\cup U(X_{\lg \rho\rg})$ and
$$
\begin{array}{l}
W(X_{\lg \rho\rg})=\{B_{(i,j)}\mid i,j\in\ZZ_p\}\,\,\,\, \text{with}\,\,\,\, B_{(i,j)}=\{(i,j,k)\mid k\in\ZZ_p\},\\
U(X_{\lg \rho\rg})=\{B_{[x,y]}\mid x,y\in\ZZ_p\} \,\,\,\, \text{with}\,\,\,\,
B_{[x,y]}=\{[x,y,z]\mid z\in\ZZ_p\}.
\end{array}$$
By the definition of $X$ we have that $X_{\lg \rho\rg}$ has edge set
$$
\begin{array}{lll}
E(X_{\lg \rho\rg})&=&\{\{B_{(i,j)},B_{[x,i+b]}\}\mid i,j,x\in\ZZ_p,b\in \Lambda\}\cup\\
&&\{\{B_{(i,j)},B_{[x,j+b]}\}\mid,i,j,x\in\ZZ_p,b\in \Lambda\},
\end{array}
$$
where $\Lambda$ is given in Definition \ref{X}.
For each $i\in\mathbb{Z}_p$, define $\mathcal{P}_i$ as the following path in the quotient graph $X_{\langle \rho \rangle}$:
$$\mathcal{P}_i: B_{[0,i+b]}\sim B_{(i,i)}\sim B_{[1,i+b]}\sim B_{(i,i+1)}\sim\cdots \sim B_{[p-1,i+b]}\sim B_{(i,i+p-1)}.$$
Then we obtain a Hamilton cycle $\mathcal{C}_{\lg \rho\rg}$ of $X_{\lg \rho\rg}$,
where
$$\mathcal{C}_{\lg \rho\rg}:\mathcal{P}_0\sim \mathcal{P}_{p-1}\sim\mathcal{P}_{p-2}\sim
\cdots\sim\mathcal{P}_{1}\sim B_{[0,b]}.$$
Consider the induced graph $X[B_{[0,b]},B_{(0,0)}]$. We have $(0,0,k)\sim [0,b,k+\frac{p\pm1}{2}]$ for any $k\in\ZZ_p$. This implies that $d(B_{[0,b]},B_{(0,0)})\geq 2$.
Then, by Proposition~\ref{pro:4}, $X$ has a Hamilton cycle.

\vskip 3mm
{\it Case 2: $s=\th^{\frac{p-1}{r}}$.}
\vskip 3mm
Now $X=X_2(p,r)$. Define the quotient graph $X_{\lg \rho\rg}$ which is induced by the automorphism $\rho$, where $\rho:(i,j,k)\mapsto (i+1,j+1,k)$
and $[x,y,z]\mapsto [x,y+1,z]$, with $i,j,k,x,y,z\in\ZZ_p$ and $\lg \rho\rg\cong\ZZ_p$.
Then we obtain a quotient graph $X_{\lg \rho\rg}$ induced by $\lg \rho\rg$, where $V(X_{\lg \rho\rg})=W(X_{\lg \rho\rg})\cup U(X_{\lg \rho\rg})$ and
$$
\begin{array}{l}
W(X_{\lg \rho\rg})=\{B_{(\delta,k)}\mid \delta,k\in\ZZ_p\}\,\,\,\, \text{with}\,\,\,\, B_{(\delta,k)}=\{(i,i+\delta,k)\mid i\in\ZZ_p\},\\
U(X_{\lg \rho\rg})=\{B_{[x,z]}\mid x,z\in\ZZ_p\} \,\,\,\, \text{with}\,\,\,\,
B_{[x,z]}=\{[x,y,z]\mid y\in\ZZ_p\}.
\end{array}$$
By the definition of $X$ we have that $X_{\lg \rho\rg}$ has edge set
$$
\begin{array}{lll}
E(X_{\lg \rho\rg})&=&\{\{B_{(\delta,k)},B_{[x,k+\frac{p+1}{2}]}\}\mid \delta,k,x\in\ZZ_p\}\cup\{\{B_{(\delta,k)},B_{[x,k+\frac{p-1}{2}]}\}\mid \delta,k,x\in\ZZ_p\}.
\end{array}
$$
Note that $\mathcal{P}_t$ is a path of the quotient graph $X_{\lg \rho\rg}$ for each $t\in\ZZ_p$, where
$$\mathcal{P}_t: B_{[0,\frac{p+1}{2}+t]}\sim B_{(0,t)}\sim B_{[1,\frac{p+1}{2}+t]}\sim B_{(1,t)}\sim\cdots\sim B_{[p-1,\frac{p+1}{2}+t]}\sim B_{(p-1,t)}.$$
Then we obtain a Hamilton cycle $\mathcal{C}_{\lg \rho\rg}$ of $X_{\lg \rho\rg}$ given by
$$\mathcal{C}_{\lg \rho\rg}:\mathcal{P}_0\sim \mathcal{P}_{p-1}\sim\mathcal{P}_{p-2}\sim
\cdots\sim\mathcal{P}_{1}\sim B_{[0,\frac{p+1}{2}]}.$$
Consider the induced graph $X[B_{[0,\frac{p+1}{2}]},B_{(0,0)}]$. We have $(0,0,0)\sim(0,sb,\frac{p+1}{2}),$ where $s=\theta^{\frac{p-1}{2r}},$ $\ZZ_p^*=\lg \theta\rg$, $b\in\Lambda\leq \ZZ_p^*$ and $|\Lambda|=r\geq 2$.
This implies that $d(B_{[0,\frac{p+1}{2}]},B_{(0,0)})\geq 2$. Then, by Proposition~\ref{pro:4}, $X$ has a Hamilton cycle.
\qed

\end{document}